\documentclass[11pt]{amsart}

\usepackage[top=1in, left=1in, right=1in, bottom=1in]{geometry}

\usepackage{amsmath,amssymb,amsthm,amsfonts}
\usepackage{paralist}
\usepackage{booktabs}
\usepackage[usenames,dvipsnames]{color}
\usepackage{comment}
\usepackage{graphicx,epsf}
\usepackage{graphics}
\usepackage{epsfig}
\usepackage{epstopdf}
\usepackage{psfrag}

\usepackage{amsmath,amsfonts,amssymb,latexsym,amsthm}
\usepackage{graphicx,epsf}
\usepackage{amsthm}
\usepackage{multicol}
\usepackage{ytableau}

\usepackage{url}
\usepackage{hyperref} 

\usepackage{tikz}
\usetikzlibrary{automata,positioning,arrows}

\usepackage[colorinlistoftodos]{todonotes}

\RequirePackage{cleveref}
\usepackage{hypcap}
\hypersetup{colorlinks=true, citecolor=darkblue, linkcolor=darkblue}
\definecolor{darkblue}{rgb}{0.0,0,0.7}
\newcommand{\darkblue}{\color{darkblue}}

\definecolor{darkred}{rgb}{0.68,0,0}

\definecolor{darkgreen}{rgb}{0,.38,0}
\newcommand{\darkgreen}{\color{darkgreen}}

\newcommand{\defn}[1]{\emph{\darkblue #1}}

\newcommand{\defng}[1]{\emph{\darkgreen #1}}

\newlength{\ml}
\newtheorem{thm}{Theorem}[section]

\newtheorem{lemma}[thm]{Lemma}

\newtheorem{cor}[thm]{Corollary}
\newtheorem{prop}[thm]{Proposition}

\newtheorem{conj}[thm]{Conjecture}

\newtheorem{op}[thm]{Open Problem}
\newtheorem{ex}[thm]{Example}

\theoremstyle{plain}

\theoremstyle{definition}

\newtheorem{rem}[thm]{Remark}

\numberwithin{equation}{section}

\def\bx{{\textbf{\textit{x}}}}
\def\by{{\textbf{\textit{y}}}}

\newcommand{\la}{\lambda}

\newcommand{\cP}{\mathcal P}

\newcommand{\nn}{\Bbb N}

\newcommand{\cc}{\Bbb C}

\newcommand{\ga}{\gamma}

\newcommand{\de}{\delta}

\newcommand{\al}{\alpha}
\newcommand{\be}{\beta}

\newcommand{\ssu}{\subset}

\def\per{\mathrm{per}}
\def\sgn{\mathrm{sgn}}

\def\.{\hskip.06cm}
\def\ts{\hskip.03cm}

\def\SSYT{ {\text {\rm SSYT}  } }

\def\GL{ {\text {\rm GL}  } }
\def\rg{\overline{g}}

\def\P{{{\rm{\textsf{P}} }}}

\def\FP{{\rm{\textsf{FP}}}}
\def\PH{{\rm{\textsf{PH}}}}

\def\SP{{{\rm{\textsf{\#P}}}}}

\newcommand{\ComCla}[1]{\textup{\textsf{#1}} }

\def\GapP{{{\rm{\textsf{GapP}} }}}
\newcommand{\CeqP}{\ComCla{C$_=$P}}
\def\FP{{{}\rm{\textsf{FP}} }}
\def\NP{{{\rm{\textsf{NP}} }}}
\def\coNP{{{\rm{\textsf{coNP}} }}}

\def\VP{{{\rm{\textsf{VP}} }}}
\def\VNP{{{\rm{\textsf{VNP}} }}}

\def\IC{{\mathbb{C}}}

\def\VPs{\mathsf{VP}_{\mathrm{ws}}}
\def\dc{\mathsf{dc}}

\def\<{\langle}
\def\>{\rangle}

\def\det{\mathrm{det}}
\def\per{\mathrm{per}}
\def\sgn{\mathrm{sgn}}
\def\GL{\mathrm{GL}}

\def\la{\lambda}

\def\id{\mathrm{id}}

\title[Structure constants]{Complexity and asymptotics of structure constants }

\author{Greta Panova}
\thanks{ \hspace{-0.25in} Department of Mathematics, University of Southern California, 
\url{https://sites.google.com/usc.edu/gpanova/} \\ To appear in \emph{Proceedings of the Open Problems in Algebraic Combinatorics} conference, May 2022. \url{https://www.samuelfhopkins.com/OPAC/opac.html}}

\begin{document}
\maketitle

\dedicatory{{\hfill {\footnotesize {``Thinking out of the black box''}}}}


\begin{abstract}

Kostka, Littlewood-Richardson, Kronecker, and plethysm coefficients are fundamental quantities in algebraic combinatorics, yet many natural questions about them stay unanswered for more than 80 years. Kronecker and plethysm coefficients lack ``nice formulas'', a notion that can be formalized using computational complexity theory. Beyond formulas and combinatorial interpretations, we can attempt to understand their asymptotic behavior in various regimes, and inequalities they could satisfy. Understanding these quantities has applications beyond combinatorics. On the one hand, the asymptotics of structure constants is closely related to understanding the [limit] behavior of vertex and tiling models in statistical mechanics. More recently, these structure constants have been involved in establishing computational complexity lower bounds and separation of complexity classes like VP vs VNP, the algebraic analogs of P vs NP in arithmetic complexity theory. Here we discuss the outstanding problems related to asymptotics, positivity, and complexity of structure constants focusing mostly on the Kronecker coefficients of the symmetric group and, less so, on the plethysm coefficients.

This expository paper is based on the talk presented at the Open Problems in Algebraic Combinatorics coneference in May 2022.

\end{abstract}

\section{Introduction}

Algebraic Combinatorics studies symmetries via their manifestations in discrete objects, connecting with areas like Representation Theory, Statistical Mechanics, Computational Complexity Theory, and Algebraic Geometry. Many of its problems arise from the need for a quantitative and explicit understanding of algebraic phenomena like group representations and decompositions into irreducible representations, dimension formulae for modules, intersection numbers in geometry, etc.

In its origin lies the representation theory of the symmetric group $S_n$ and general linear group $GL_N$. Their irreducible modules are indexed by integer partitions, and their bases can be described via standard and semi-standard Young tableaux (SYTs and SSYTs). As a next step, it is natural to understand how tensor products of irreducible representations decompose into irreducible components. In the case of $GL_N$ the multiplicities of isotypic components are the Littlewood-Richardson (LR) coefficients. For $S_n$ the multiplicities are the Kronecker coefficients. The composition of irreducible $GL$ modules decomposes into irreducibles with multiplicities given by the plethysm coefficients. 

While no compact formula for the LR coefficients exists, they can be understood via their combinatorial interpretation, the LR rule, which gives them in the form of certain skew SSYTs. The LR rule was formulated in 1934, but it took 40 years to prove formally. Inspired by the LR coefficients, Murnaghan defined the Kronecker coefficients in 1938 and observed that computing them is a highly nontrivial task. In late 70s and 80s the research turned towards obtaining positive formulas. In 2000 Stanley stated  as a major open problem in Algebraic Combinatorics the problem of finding a nonnegative combinatorial interpretation for the Kronecker coefficients, as well as the plethysm coefficients. However, the progress has been minimal. 

Besides presenting fundamental natural questions needing answers, the Kronecker and plethysm coefficients also play an important role in algebraic complexity theory, specifically Geometric Complexity Theory (GCT). Algebraic complexity theory studies the complexity of computing formal polynomials using arithmetic circuits, and some of its main problems concern the separation of complexity classes $\VP \neq \VNP$, the algebraic analog of the $\P \neq \NP$ Millennium problem. GCT approaches these problems by studying the symmetries and geometric properties of the universal polynomials (usually determinant and permanent), passing onto the coordinate rings of $GL$ orbit closures of these polynomials. Complexity lower bounds and separation of classes then reduces to understanding the $GL$-irreducible representations of these rings and comparing their multiplicities. Such multiplicities are closely related to LR, Kronecker, and plethysm coefficients, and separating complexity classes can be achieved by proving certain inequalities for such multiplicities. 

Another application of Algebraic Combinatorics is in Statistical Mechanics, Asymptotic Representation Theory and Random Matrix Theory. Some of the early connections date to the remarkable story of the longest increasing subsequence of permutations, which has the same [Tracy-Widom] distribution as the maximal eigenvalue of a random Hermitian matrix and the height of an SYT from the Plancherel measure.  Symmetric polynomials, and particularly Schur functions, appear extensively in the study of random lozenge tilings, exclusion processes and vertex models. Plane partitions themselves are lozenge tilings with certain boundaries. Conversely, probabilistic and asymptotic methods have been used to understand algebro-combinatorial quantities. In the lack of compact formulas and combinatorial interpretations, understanding the growth behavior of structure constants is a natural next step. Beyond the interaction with Statistical Mechanics, understanding their growth is closely related to their application in GCT mentioned above. 

\subsection*{Paper structure}
Here we will briefly define the relevant structure constants, give some background, and phrase the open problems of three  separate aspects: positivity, complexity, asymptotics.   
For more background on Complexity Theory and its specific connections with Algebraic Combinatorics we refer to the companion paper~\cite{Pan23} and references therein.

This paper presents the professional view of the author and is centered around problems on which she has worked on. It has no aspirations to be a comprehensive survey on the topics.  

\subsection*{Acknowledgements}
The author is grateful to her collaborators throughout the years and in particular Christian Ikenmeyer and Igor Pak for the extensive work on the subject of this paper. We would like to dedicate this to the memory of Christine Bessenrodt, whose deep and insightful work in the area has led a lot of the progress and with whom the author had the honor to collaborate briefly in 2020. 
We are grateful to the organizers of ``Open Problems in Algebraic Combinatorics''  for the inspiring conference and the opportunity to present these problems.  

The author has been partially supported by the NSF.

\section{Basic objects and notions}\label{s:ac}

We will assume the reader is familiar with basic notions in Algebraic Combinatorics like integer partitions, Young tableaux and symmetric functions. For a concise recap of those see~\cite{Pan23}. For details on the combinatorial sides see~\cite{S1,Mac} and for the representation theoretic aspects see~\cite{Sag,Ful}.

\subsection{Partitions and Young tableaux} \label{ss:basics}
We use standard notation from~\cite{Mac} and~\cite[$\S$7]{S1}
throughout the paper.

Let $\la=(\la_1,\la_2,\ldots,\la_\ell)$ be a \defn{partition}
of size $n:=|\la|=\la_1+\la_2+\ldots+\la_\ell$, where
$\la_1\ge \la_2 \ge \ldots \ge\la_\ell\ge 1$.  We write
$\la\vdash n$ for this partition, and $\cP=\{\la\}$ for
the set of all partitions.   The length of~$\la$
is denoted $\ell(\la):=\ell$. Denote by $p(n)$ the number
of partitions $\la\vdash n$. Let $\la+\mu$ denote the partition
$(\la_1+\mu_1,\la_2+\mu_2,\ldots)$

Special partitions include the \emph{rectangular shape} \ts
$(a^b) = (a,\ldots,a)$, $b$ times, the \emph{hooks shape}
$(k,1^{n-k})$, the \emph{two-row shape} $(n-k,k)$, and the
\emph{staircase shape} \ts $\rho_\ell = (\ell,\ell-1,\ldots,1)$.

\ytableausetup{boxsize=2ex}

A \defn{Young diagram} of
\emph{shape}~$\la$ is an arrangement of squares
$(i,j)\ssu \nn^2$ with $1\le i\le \ell(\la)$
and $1\le j\le \la_i$.  Let $\la\vdash n$.
A \defn{semistandard Young tableau} $T$ of
\emph{shape}~$\la$ and \emph{weight}~$\al$ is an
arrangement of $\al_k$ many integers~$k$ in squares
of~$\la$, which weakly increase along rows and strictly
increase down columns, i.e. $T(i,j) \leq T(i,j+1)$ and $T(i,j) \leq T(i+1,j)$. For example,
\ytableausetup{centertableaux} 
\ytableaushort{11244,2235,45}\, is an SSYT of shape $\la=(5,4,2)$ and type $\al=(2,3,1,3,2)$.
Denote by $\SSYT(\la,\al)$ the
set of such tableaux, and \ts
$K(\la,\al) = \bigl|\SSYT(\la,\al)\bigr|$ \ts the
\defn{Kostka number}. A \defn{standard Young tableau} (SYT) of shape $\la$ is an SSYT of type $(1^n)$, and we have $f^\la := K_{\la,1^n}$, which can be computed by the hook-length formula (HLF) of~\cite{FRT}.

\subsection{Representations of $S_n$ and $GL_N$}
The irreducible representations of the \emph{symmetric group} $S_n$ are the \defn{Specht modules} $\mathbb{S}_\la$ and are indexed by partitions $\la \vdash n$. A basis for $\mathbb{S}_\la$ can be indexed by the SYTs. In particular
$$\dim \mathbb{S}_\la = f^\la.$$
We have that $\mathbb{S}_{(n)}$ is the trivial representation assigning to every $w$ the value $1$ and $\mathbb{S}_{1^n}$ is the sign representation. 

The \defn{character $\chi^\la(w)$} of $\mathbb{S}_{\la}$ can be computed via the \emph{Murnaghan-Nakayama} rule. Let $w$ have type $\al$, i.e. it decomposes into cycles of lengths $\al_1,\al_2,\ldots,\al_k$. Then
$$\chi^\la(w)=\chi^\la(\al) = \sum_{T \in MN(\la;\al)} (-1)^{ht(T)},$$
where $MN$ is the set of rim-hook tableaux of shape $\la$ and type $\al$, so that the entries are weakly increasing along rows and down columns, and all entries equal to $i$ form a rim-hook shape (i.e. connected, no $2\times 2$ boxes) of length $\al_i$. The height of each rim-hook is one less than the number of rows it spans, and $ht(T)$ is the sum of all these heights. For examples
$$\ytableaushort{112333,12234,22334}$$
is a Murnaghan-Nakayama tableau of shape $(6,5,5)$, type $(3,5,6,2)$ and has height $ht(T) = 1 + 2 + 2 +1=6$.


The irreducible polynomial representations of  $GL_N(\mathbb{C})$ are the \emph{Weyl modules} $V_\la$ and are indexed by all partitions with $\ell(\la) \leq N$. Their characters are exactly the Schur functions $s_\la(x_1,\ldots,x_N)$, where $x_1,\ldots,x_N$ are the eigenvalues of $g \in GL_N(\mathbb{C})$.


\subsection{Symmetric functions}\label{ss:sym_func}
Let $\Lambda[\bx]$ be the ring of \defn{symmetric functions} $f(x_1,x_2,\ldots)$,
where the symmetry means that $f(\bx)=f(\bx_\sigma)$ for any permutation $\sigma$
of the variables, and $f$ is a formal power series.
The ring $\Lambda_n$ of homogeneous symmetric functions of degree $n$ has several important bases: the \defng{homogenous symmetric functions} \. $h_\la$\.,
\defng{elementary symmetric functions} \. $e_\la$,
\defng{monomial symmetric functions} \. $m_\la$ \.,
\defng{power sum symmetric functions} \. $p_{\la}$   and
\defng{Schur functions} \. $s_\la$. 
The Schur functions can be defined as the generating functions for SSYTs of shape $\la$
$$s_\la = \sum_{\mu \vdash n} K_{\la\mu} m_\mu,$$
but also using \emph{Weyl's determinantal} formula
$$s_\la(x_1,\ldots,x_\ell) = \frac{ \det [x_{i}^{\la_j+\ell-j}]_{i,j=1}^\ell}{\prod_{i<j} (x_i-x_j)}$$
and the \emph{Jacobi-Trudi} identity
$$s_\la = \det [h_{\la_i -i+j}]_{i,j=1}^{\ell(\la)}.$$

The \emph{Cauchy identity}
$$\sum_\la s_\la(\mathbf{x})s_\la(\mathbf{y}) = \prod_{i,j} \frac{1}{1-x_iy_j}$$
has a remarkable combinatorial proof given by the RSK correspondence between pairs of same shape SSYTs and matrices with nonnegative integer entries.

\subsection{Computational Complexity}\label{ss:cc}
We refer to \cite{Aa,Wig,BCS} for details on Computational Complexity classes, and to~\cite{Pak22,Pan23} and references therein for the connections with Algebraic Combinatorics.

A \defn{decision problem} is a computational problem, for which the output is Yes or No. There are two major complexity classes $\P$ and~$\NP$, subject of the $\P$ vs $\NP$ Millennium problem. $\P$ is the class of decision problems, where given any input of size $n$ (number of bits required to encode it),  such that the answer can be obtained in \emph{polynomial time} denoted by $poly(n)$, i.e. there is a fixed $k$ and an algorithm taking $O(n^k)$ many steps (elementary operations). $\NP$ is the class of decision problems, where if the answer is Yes, then it can be verified in polynomial time, i.e. there is a poly-time computable witness.  Naturally, $\P \subset \NP$ and it is widely believed that $\P \neq \NP$.
The classes $\FP$ and $\SP$ are the counting analogues of $\P$ and $\NP$. A \emph{counting problem} is in $\SP$ if it is the number of accepting paths of an $\NP$ Turing machine. In practice, these are counting problems where the answers are exponentially large sums of $0$-$1$ valued functions $M$, each of which can be computed in $O(n^k)$ time for a fixed $k$:
$$\sum_{b \in \{0,1\}^{m}} M(b),$$ 
where $m=poly(n)$.  We set the classes $\GapP = \{ f-g |\, f,g \in \SP\}$ and $\GapP_{\geq 0} = \{ f-g | \, f,g \in \SP \text{ and } f-g \geq 0\}$.

\smallskip

The algebraic complexity classes $\VP$ and $\VNP$ were introduced by Valiant~\cite{V1}, as the algebraic analogues of $\P$ and $\NP$ (we refer to~\cite{Bur00} for formal definitions and properties). They concern the computation of polynomials $f \in \mathbb{F}[x_1,\ldots,x_N]$ of degree $poly(N)$ using arithmetic circuits where the inputs are constants of arbitrary size from the field $\mathbb{F}$ and the formal variables $x_1,\ldots,x_n$, the gates are $+, \times$, and the output is $f$. Then $\VP$ is the class of polynomials for which there is a fixed $k$ and a circuit of size $O(n^k)$ computing the polynomial $f$. $\VNP$ is the class of polynomials $f$, such that there exist a $k$ and $m = O(N^k)$ and a polynomial $g \in \VP$ in $N+m$ variables, such that
$$f(x_1,\ldots,x_N) = \sum_{b \in \{0,1\}^m} g(x_1,\ldots,x_N,b_1,\ldots,b_m).$$

The main conjecture is that $\VP \neq \VNP$ and is closely related to $\P \neq \NP$, see e.g.~\cite{Bur0}.  As Valiant showed, for every polynomial $f \in \mathbb{C}[x_1,\ldots,x_n]$  there exists a $K$ and a $K \times K$ matrix $A := A_0 + \sum_{i=1}^N A_i x_i$ with $A_j \in \mathbb{C}^{K \times K}$  such that $\det A= f$.  The smallest such $K$ is the \defn{determinantal complexity} of $f$ denoted $\dc(f)$ and it is finite for every $f$.
The $\VPs$-universal polynomial is the determinant, in the sense that $f \in \VPs$ iff $\dc(f) = poly(n)$, where $poly(n)$ denotes any fixed degree polynomial in $n$. Here $\VPs$ is the class of polynomials with poly-sized Algebraic Branching Programs (a model between ``formula'' and ``circuit'''). The class $\VPs$ has been the main subject of study and is often interchangably used for $\VP$. 
The classical universal $\VNP$-complete polynomial is the permanent
{\small
$$\per_m [X_{ij}]_{i,j=1}^m = \sum_{\sigma \in S_m} \prod_{i=1}^m X_{i\sigma(i)}.$$
}
Thus, separating $\VPs$ from $\VNP$ can be done by showing that $\dc(per_m) \gg poly(m)$.

The \defn{Geometric Complexity Theory} (GCT) program of Mulmuley and Sohoni~\cite{MS1,MS2} compares the determinant and permanent by comparing their orbit closures under the action of $GL_{n^2}$ on the variables $(X_{11},\ldots,X_{nn})$. If $\dc(\per_m) \leq n$, then $\overline{GL_{n^2} \per_m} \subset \overline{ GL_{n^2}{\det}_n}$, and passing to the coordinate ring, this induces a $GL_{n^2}$-equivariant surjection
\begin{align}\label{eq:inclusion}
   \mathbb{C}[\overline{GL_{n^2} {\det}_n}]  \twoheadrightarrow \mathbb{C}[\overline{GL_{n^2} X_{11}^{n-m} \per_m}]
\end{align} To show a lower bound for $\dc(\per_m)$, i.e. $\dc(\per_m)>n$, then we need to show that the above surjection does not happen. Decomposing these modules into irreducibles for any given degree $d$, we have
\begin{equation}\label{eq:coord_rings_mult}
\IC[\overline{\GL_{n^2} {\det}_n }]_d \simeq \bigoplus_{\la } V_\la^{\oplus \delta_{\lambda,d,n}}, \quad\text{and}\quad
\IC[\overline{\GL_{n^2} \per_m^n}]_d \simeq \bigoplus_{\la} V_\la^{\oplus \gamma_{\lambda,d,n,m}}.
\end{equation}

 Thus, to show that~\ref{eq:inclusion} is not possible,  we can find an irreducible $GL_{n^2}$ module $V_\la$, such that $ \gamma_{\lambda,d,n,m} > \delta_{\la,d,n}$ for some $d$. Such $\la$'s are called \defn{multiplicity obstructions} and $\ga$ and $\delta$ are closely related to plethysm and Kronecker coefficients. When $\delta_{\la,d,n}=0$ also, then $\la$ is an \defn{occurrence obstruction}. It was shown in~\cite{BIP} that there are no occurrence obstructions for $n>m^{25}$, so we need to find multiplicity obstructions.
 
 This is one of the main applications of Kronecker and plethysm coefficients and motivates our work. In particular,  we have 
$$\de_{\la,d,n} \leq g(\la,n^d, n^d) \qquad \text{and} \qquad \ga_{\la,d,n,m} \leq a_\la(d[n]).$$
See \S\ref{s:mult} for the definitions
 and~\cite{Pan23} for a brief overview of these connections.

\section{Multiplicities and structure constants}\label{s:mult}

Once the irreducible representations have been sufficiently understood, it is natural to consider what other representations can be formed by them and how such representations decompose into irreducibles. 

In the case of $GL_N(\mathbb{C})$ these coefficients are
the  \defn{Littlewood--Richardson} \defn{coefficients} (LR) $c^{\la}_{\mu\nu}$ defined as the multiplicity of $V_\la$ in $V_\mu \otimes V_\nu$, so
$$V_\mu \otimes V_\nu = \bigoplus_\la V_\la ^{\oplus c^{\la}_{\mu\nu}}.$$
Via their characters, they can be equivalently defined as

$$s_\mu s_\nu = \sum_\la c^{\la}_{\mu\nu} s_{\la}, \quad \text{ and } \quad s_{\la/\mu} = \sum_\nu c^{\la}_{\mu\nu}s_\nu.$$

While no nice product formula exists for their computation, they have a \emph{combinatorial interpretation}, the so called \emph{Littlewood-Richardson rule}, see~\cite[App. 1]{S1}.

\begin{thm}[Littlewood-Richardson rule]
The Littlewood-Richardson coefficient $c^{\la}_{\mu\nu}$ is equal to the number of skew SSYT $T$ of shape $\la/\mu$ and type $\nu$, whose reading word is a ballot sequence.
\end{thm} 

A ballot sequence means that in every prefix of the sequence the number of $1$'s$\geq$ the number of $2$'s  $\geq$ number of $3$'s etc.
For example, we have $c^{(6,4,3)}_{(3,1),(4,3,2)}=2$ as there are two LR tableaux (SSYT of shape $(6,4,3)/(3,1)$ and type $(4,3,2)$ whose reading words are ballot sequences)
\ytableaushort{\none\none\none 111,\none 122,233}\, with reading word $111221332$ and \ytableaushort{\none\none\none 111,\none 222,133}\,   with reading word 111222331.

The \defn{Kronecker coefficients} $g(\la,\mu,\nu)$ of the symmetric group are the corresponding structure constants for the ring of $S_n$- irreducibles. Namely, $S_n$ acts diagonally on the tensor product of two Specht modules and the corresponding module factors into irreducibles with multiplicities $g(\la,\mu,\nu)$
$$\mathbb{S}_\la \otimes \mathbb{S}_\mu = \bigoplus_\nu \mathbb{S}_\nu^{\oplus g(\la,\mu,\nu)}.$$
In terms of characters we can write them as 
\begin{equation}\label{eq:char_kron}
g(\la,\mu,\nu) = \langle \chi^\la \chi^\mu,\chi^\nu\rangle = \frac{1}{n!} \sum_{w\in S_n} \chi^\la(w) \chi^\mu(w) \chi^\nu(w).
\end{equation}
The last formula shows that they are symmetric upon interchanging the underlying partitions $g(\la,\mu,\nu)= g(\mu,\nu,\la) = \cdots$ which motivates us to use such symmetric notation. 

The \defn{Kronecker product} $*$ on $\Lambda$ is defined on the Schur basis by 
$$s_\la * s_\mu = \sum_\nu g(\la,\mu,\nu) s_\nu,$$
and extended by linearity. 

Inspired by the Littlewood-Richardson coefficients, Murnaghan defined them in 1938. In fact, he showed that, see~\cite{Mur56},
\begin{thm}[Murnaghan]\label{thm:mur}
For every $\la,\mu,\nu$, such that $|\la|=|\mu|+|\nu|$ we have that 
$$c^\la_{\mu\nu} = g( (n - |\la|,\la), (n-|\mu|,\mu), (n-|\nu|,\nu)),$$
for sufficiently large $n$.
\end{thm}
In particular, one can see that $n=2|\la|+1$ would work. Note that even when the sizes of the partitions do not add up, the RHS stabilizes as $n \to \infty$ to the \emph{reduced Kronecker coefficients}, see~\ref{s:reduced}. 

Thanks to the Schur-Weyl duality, they can also be interpreted via Schur functions as
$$s_\la[\bx \by] = \sum_{\mu,\nu} g(\la,\mu,\nu) s_\mu(\bx) s_\nu(\by),$$
where $\bx\by = (x_1y_1,x_1y_2,\ldots,x_2y_1,\ldots)$ is the vector of all pairwise products. In terms of $GL$ representations they give us the dimension of the invariant space
$$g(\la,\mu,\nu) = \dim ( V_\mu \otimes V_\nu \otimes V^*_\la)^{GL_{N}\times GL_M},$$
where $V_\mu$ is considered as a $GL_N$-module, and $V_\nu$ as a $GL_M$-module and $V_\la^*$ is a dual $GL_{NM}$-module. In other words, they are dimensions of highest weight vector spaces, and using  this interpretation for them as dimensions of highest weight spaces one can show the following, see~\cite{CHM}
\begin{thm}[Semigroup property]\label{thm:semigroup}
Let $(\al^1,\be^1,\ga^1)$ and $(\al^2,\be^2,\ga^2)$ be two partition triples, such that $|\al^1|=|\be^1|=|\ga^1|$ and $|\al^2|=|\be^2|=|\ga^2|$. Suppose that $g(\al^1,\be^1,\ga^1)>0$ and $g(\al^2,\be^2,\ga^2)>0$. Then
$$g(\al^1+\al^2,\be^1+\be^2,\ga^1+\ga^2) \geq \max\{ g(\al^1,\be^1,\ga^1), g(\al^2,\be^2,\ga^2)\}.$$
\end{thm}

From the original $S_n$ characters definition we see that
$$g(\la,\mu,\nu) = g(\la',\mu',\nu)$$
since $\chi^{\la'}(w) = \chi^{\la}(w)\sgn(w)$. In particular, we have $g(\la,\mu,(n))=\delta_{\la,\mu}$ for all $\la$ and $g(\la,\mu,1^n)=\delta_{\la,\mu'}$. 

\begin{ex}
By the above observation we have that  
$$h_k[\bx\by] =s_k[\bx\by] = \sum_{\la \vdash k} s_\la(\bx) s_\la(\by).$$
Using the Jacobi-Trudi identity we can write
\begin{align*}
s_{2,1}[\bx\by] = h_2[\bx\by]h_1[\bx\by] - h_3[\bx\by] = \left(s_2(\bx)s_2(\by) +s_{1,1}(\bx)s_{1,1}(\by)\right)s_1(\bx)s_1(\by) \\
- s_3(\bx)s_3(\by) -s_{2,1}(\bx)s_{2,1}(\by) - s_{1,1,1}(\bx)s_{1,1,1}(\by)\\
= s_{2,1}(\bx)s_{2,1}(\by) + s_{2,1}(\bx)s_3(\by) +s_3(\bx)s_{2,}(\by) +s_{1,1,1}(\bx)s_{2,1}(\by) + s_{2,1}(\bx)s_{1,1,1}(\by).
\end{align*}
So we see that $g((2,1),(2,1),(2,1))=1$. 
\end{ex}

The \defn{plethysm coefficients} $a^\la_{\mu,\nu}$ are multiplicities of an irreducible $\GL$ representation in the composition of two $\GL$ representations. Namely, let $\rho^\mu: GL_N \to GL_M$ be one irreducible, and $\rho^\nu: GL_M \to GL_K$ be another. Then $\rho^\nu \circ \rho^\mu: GL_N \to GL_K$ is another representation of $GL_N$ which has a character $s_\nu[s_\mu]$, which decomposes into irreducibles as
$$s_\nu[s_\mu] = \sum_\la a^\la_{\mu,\nu}s_\la.$$
Here the \emph{plethystic} notation $f[g]$ means the evaluation of $f$ over the monomials of $g$ as variables: if $g = \bx^{\al^1} + \bx^{\al^2} +\cdots$, then $f[g] = f(\bx^{\al^1},\bx^{\al^2},\ldots)$. 

\begin{ex}
We have that
\begin{align*}
s_{(2)}[s_{(1^2)}] = h_3[e_2] = h_2(x_1x_2,x_1x_3,\ldots)=x_1^2x_2^2 + x_1^2x_2x_3 + 3x_1x_2x_3x_4+\cdots \\
= s_{2,2}(x_1,x_2,x_3,\ldots) + s_{1,1,1,1}(x_1,x_2,x_3,\ldots),  
\end{align*}
so $a^{(2,2)}_{(2),(1,1)} =1$ and $a^{(3,1)}_{(2),(1,1)}=0$.
\end{ex}

We will be particularly interested when $\mu = (d)$ or $(1^d)$ which are the $d$th symmetric power $Sym^d$ and the $d$th wedge power $\Lambda^d$, and $\nu = (n)$. We denote this plethysm coefficient by $a_\la(d[n]) := a^\la_{(d),(n)}$ and
$$h_d[h_n] = \sum_\la a_\la(d[n]) s_\la.$$
These coefficients are central to the GCT application and deserve special attention.

\section{Combinatorial interpretations}\label{s:combin}

``Combinatorial interpretation'', albeit never formally defined in the literature, is assumed to mean a family of nice discrete objects whose cardinality gives the desired quantity. In practice,  for quantities with combinatorial interpretation, computing them is in \SP. On the other hand, being in \SP\ does not necessarily give a ``nice'' combinatorial interpretation, but is certainly a starting point. 
All problems considered here are in $\GapP$, which makes $\SP$ the natural class to consider in relation to ``manifestly nonnegative'' formulas. 
If such a problem is not in $\SP$, then there could not be a reasonable nonnegative combinatorial interpretation within the given context. This paradigm is discussed at length in~\cite{IP22,Pak22,Pan23}. 

Following the discovery of the Littlewood-Richardson rule in 1934, Murnaghan~\cite{Mur38} defined the Kronecker coefficients of $S_n$ and computed various special cases, but quickly observed that this is no easy task. Interest in nonnegative combinatorial interpretation of these coefficients reemerged in the 80s, see~\cite{Las79,GR85}, and was stated as Problem 10 in Stanley's list\footnote{See this for the original list and updates on the problems \url{https://mathoverflow.net/questions/349406/}} of ``Open Problems in Algebraic Combinatorics''~\cite{Sta00}. 

\begin{op}[Stanley]\label{op:kron}
Find a combinatorial interpretation of $g(\la,\mu,\nu)$, thereby combinatorially proving that they are nonnegative. 
\end{op}

Over the years, there has been very little progress on the question. In 1989 Remmel determined $g(\la,\mu,\nu)$ when two of the partitions are hooks. In 1994 Remmel and Whitehead~\cite{RW} determined $g(\la,\mu,\nu)$ when two of the partitions are two-rows, i.e. $\ell(\la), \ell(\mu)\leq 2$. This case was subsequently studied also in~\cite{BMS}. In 2006 Ballantine and Orellana~\cite{BO} determined a rule for $g(\la,\mu,\nu)$ when one partition is a two-row, e.g. $\mu=(n-k,k)$, and the first row of one of the others is large, namely $\la_1 \geq 2k-1$. The most general rule was determined by Blasiak~\cite{B13} in 2012 when one partition is a hook, and this was later simplified by Blasiak and Liu~\cite{BL}, and further in~\cite{Liu17}. Other very special cases have been computed  for multiplicity free Kronecker products by Bessenrodt-Bowman~\cite{BB17}, Kroneckers corresponding to pyramids by Ikenmeyer-Mulmuley-Walter~\cite{IMW17},  $g(m^k,m^k,(mk-n,n))$ as counting labeled trees by Pak-Panova see Remark~\ref{rem:kron_trees}, near rectangular partitions by Tewari in~\cite{T15}, etc.

\medskip

Similar questions pertain to the plethysm coefficients. The following problem is number 9 in Stanley's list~\cite{Sta00}.

\begin{op}[Stanley]
Find a combinatorial interpretation for the plethysm coefficients $a^\la_{\mu,\nu}$.
\end{op}

A detailed survey on the partial results and methods can be found in~\cite{COSSZ}. 

Even the simple case for $a_\la(d[n]) = \< s_\la , h_d[h_n] \>$ is not known in general unless $n=2$. 
While, there is no direct connection between the Kronecker and plethysm coefficients, the case when $\ell(\la)=2$ coincide, see~\cite{PP14}, and actually give a combinatorial interpretation as explained below.
\begin{prop}\label{prop:kron_part}
We have that $g(\la,n^d, n^d) = a_\la(d[n])=p_{k}(n,d) - p_{k-1}(n,d)$ for $\la=(nd-k,k) \vdash nd$. Here $p_r(a,b) = \# \{ \mu \vdash r: \mu_1\leq a, \ell(\mu)\leq b\}$ are the partitions of $r$ which fit inside the $a \times b$ rectangle. 
\end{prop}

\begin{rem}\label{rem:kron_trees}
In particular, these are the coefficients in the $q$-binomials
 $$\sum_r p_r(a,b) q^r = \binom{a+b}{a}_q := \prod_{i=1}^a \frac{ (1-q^{i+b})}{1-q^i}.$$ 
 The combinatorial proof by Kathy O'Hara of the unimodality of the $p_r(a,b)$ is what gives the cominatorial interpration for the particular Kronecker and plethysm coefficients as counting certain labeled trees, see~\cite[p. 9]{Pporto}. 
\end{rem}

\bigskip

To formalize these problems, we formulate them within the Computational Complexity framework. To set the stage with an example, we start with the Littlewood-Richardson coefficients.

\medskip

\textsf{ComputeLR}: 

\texttt{Input:} $\la,\mu, \nu$ 

\texttt{Output:} Value of $c^{\la}_{\mu\nu}$.

\medskip

Here and below, the input can be encoded in \emph{unary} or \emph{binary}, and this could affect the complexity. By unary we mean that the input for each partition is written as $(\underbrace{11\ldots1}_{\la_1},\underbrace{1\ldots1}_{\la_2},\ldots)$, so the input size is $O(|\la|)=O(n)$. In binary, we encode each part $\la_i$ in base 2, so the input size is $\log_2(\la_1) + \log_2(\la_2)+\cdots =O( \ell(\la) \log_2 (\la_1))$.

By the Littlewood-Richardson rule we have that when the input is unary, \textsf{ComputeLR}$\in \SP$, as our witnesses are the skew $\la/\mu$ tableaux, which are at most $n^n = O(2^{n^2})$-many, and each can be checked to be an LR tableaux in time $O(n^2)$. The LR conditions can be translated into $O(\ell(\la)^2)$ many linear inequalities, which comprise the LR polytope. This description gives us that if the input is in binary we still have  \textsf{ComputeLR}$\in \SP$.

Using a reduction to \textsf{Knapsack} for binary input, it is not difficult to see that
\begin{thm}[Narayanan \cite{Nar}]
When the input $\la,\mu, \nu$ is in binary, \textsf{ComputeLR} is \SP-complete.
\end{thm}

When the input is in unary, we do not know whether the problem is still that hard.

\begin{conj}[Pak-Panova 2020]
When the input is in unary we have that \textsf{ComputeLR} is \SP-complete.
\end{conj}

\medskip

None of the above has been achieved for the Kronecker and plethysm coefficients, however, due to the lack of any positive combinatorial formula. 
Mulmuley also conjectured that computing the Kronecker coefficients would be in \SP,  mimicking the Littlewood-Richardson coefficients.  
\smallskip

\textsf{ComputeKron}: 

\texttt{Input:} $\la,\mu, \nu$ 

\texttt{Output:} Value of $g(\la,\mu,\nu)$.

\smallskip

\begin{conj}[Pak]
 \textsc{ComputeKron} is not in $\SP$ under reasonable complexity theoretic assumptions like $\PH$ not collapsing. 
\end{conj}

If the above is proven, that would make any solution to Open Problem~\ref{op:kron} as unlikely as the polynomial hierarchy collapsing. Any reasonable combinatorial interpretation such as counting certain objects would show that the problem is in \SP, since the objects would likely be verifiable in polynomial time.  The problem stands for both unary and binary input.

Note that \textsf{ComputeKron}$\, \in\GapP$ (\cite{BI08}) as it is easy to write an alternating sum for its computation, for example using contingency arrays, see~\cite{PP17}. 

The author's experience with Kronecker coefficients  suggests that some particular families would be as hard as the general problem.  

\begin{conj}[Panova]
 \textsf{ComputeKron} is  in \SP\  when $\ell(\la)=2$ if and only if \textsf{ComputeKron} is  in \SP\ . Likewise, \textsf{ComputeKron} is  in \SP\  for $\mu=\nu =(n^d)$ and $\la \vdash nd$ as the input if and only if \textsf{ComputeKron} is  in \SP\ .
\end{conj}

Of course, the backward direction is clear, if the general problem is in \SP, then also the special subproblems are in \SP. 


\medskip

The plethysm coefficients can also be studied from this perspective. 

\smallskip
%
%
%

\textsf{ComputePleth}: 

\texttt{Input:} $\la,\mu, \nu$

\texttt{Output:} Value of $a^\la_{\mu,\nu}$.


\begin{op}
Determine whether \textsf{ComputePleth}$\in \SP$. 
\end{op}

Using symmetric function identities, it is not hard to find an alternating formula for the plethysms and show that they are also in \GapP\!, see~\cite{FI20}. They also show that deciding positivity is \NP\!-hard and thus the computational problem is \SP-hard.  \textsf{ComputePleth} may not in be in $\SP$ in the general case, but also when $\mu,\nu$ are single row partitions. The coefficient  $a_\la(d[n])$   in this case has special significance in GCT, see~\cite{Pan23}.

\bigskip


Underlying all the representation theoretic multiplicities mentioned above are the characters of the symmetric group. For example, equation~\eqref{eq:char_kron} expresses the Kronecker coefficients via characters, and the other structure constants can also be expressed in similar ways. This motivates us to understand them also computationally. As we shall see, they also provide a proof of concept for the idea that certain positive integral algebro-combinatorial quantities may indeed not have a combinatorial interpretation. 

The characters satisfy some particularly nice identities coming from the orthogonality of the rows and columns of the character table in $S_n$.  We have that 
\begin{equation}\label{eq:char_sum}
\sum_{\la \vdash n} \chi^\la(w)^2 = \prod_i i^{c_i} c_i! \, ,
\end{equation}
where $c_i =$ number of cycles of length $i$ in $w \in S_n$. When $w=\id$, we have that $\chi^\la(\id) = f^\la$, the number of SYTs and the identity is proven using the beautiful RSK bijection. The first step in this proof is to identify $(f^\la)^2$ as the number of pairs of SYTs of the same shape.

Could anything like that be done for equation~\eqref{eq:char_sum}? The first step would be to understand what objects $\chi^\la(w)^2$ counts, does it have any positive combinatorial interpretation? We formulate it again using the CC paradigm as

\smallskip
\textsf{ComputeCharSq}: 

\texttt{Input}: $\la,\al \vdash n$, unary. 

\texttt{Output}: the integer $\chi^\la(\al)^2$.

\smallskip

\begin{thm}[\cite{IPP22}]
\textsf{ComputeCharSq}$\not \in \SP$ unless $PH=\Sigma^{\P}_2$. 
\end{thm}

The last condition says ``polynomial hierarchy collapses to the second level'', which is almost as unlikely as \P=\NP. The widely believe complexity theoretic assumptions are that $\P \neq \NP$ and that $\PH$ does not collapse. This problem shows that some natural questions in Algebraic Combinatorics are indeed not in \SP. 

The proof of this Theorem follows from another complexity theoretic result we show, namely that deciding if $\chi^\la(\al)\neq  0$ is $\PH$-hard as well, see end of \S\ref{s:positivity}.

\section{Positivity problems}\label{s:positivity}

Motivated by other developments, further questions about the Kronecker coefficients have appeared. Following their work in~\cite{HSTZ} on the square of the Steinberg character for finite groups of Lie type, the authors conjectured that a similar irreducible character should exist for the symmetric group and formulated the following. 
\begin{conj}[Tensor square conjecture]\label{conj:tensor_square}
For every $n\geq 9$ there exists a symmetric partition $\la \vdash n$, such that $\mathbb{S}_\la \otimes \mathbb{S}_\la$ contains every irreducible $S_n$ module. In other words $g(\la,\la,\mu)>0$ for all $\mu\vdash n$.
\end{conj}

Following this, Saxl conjectured a specific partition, namely the staircase, which would satisfy that.
\begin{conj}[Saxl, see~\cite{PPV}]\label{conj:saxl}
Let $\delta_k=(k,k-1,\ldots,1)$ be the staircase partition. Then $g(\delta_k,\delta_k,\mu)>0$ for all $k$ and $\mu \vdash \binom{k+1}{2}$.
\end{conj} 

This conjecture has received a lot of attention partially due to its concreteness.

Positivity results were proved using a combination of three methods -- the semigroup property constructing recursively positive triples from building blocks, explicit heighest weight constructions using the techniques in \cite{Ful00}, and an unusual comparison with characters, which was originally stated in by Bessenrodt and Behns~\cite{BB} to show that  $g(\la,\la,\la)>0$ for $\la=\la'$, later generalized in~\cite{PPV}, and in its final form in~\cite{PPq}.

\begin{thm}[\cite{PPq}]\label{prop:char_kron}
Let $\la, \mu \vdash n$ and $\la=\la'$. Let $\hat{\la} = (2\la_1-1,2\la_2-3,3\la_3-5,\ldots)$ be the principal hooks of $\la$. Then 
$$g(\la,\la,\mu)\geq |\chi^\mu(\hat{\la})|.$$
\end{thm}

The technique of explicit highest weight vectors construction led to the following more general result.

\begin{thm}[\cite{I15}]\label{thm:rho_mu}
Let $\mu \vdash \binom{k+1}{2}$ be a partition comparable to $\rho_k$ in the dominance order, i.e. $\mu \succ \rho_k$ or $\rho_k \succ \mu$. Then $g(\rho_k,\rho_k,\mu)>0$. 
\end{thm}

Using this result and generalizations, and clever combinations of the semigroup property~\ref{thm:semigroup},  Luo and Sellke~\cite{LS} proved that $\chi^{\rho_k}\otimes \chi^{\rho_k} \otimes \chi^{\rho_k} \otimes \chi^{\rho_k}$ contains every $\chi^\mu$, and that $\chi^{\rho_k}\otimes \chi^{\rho_k}$ contains ``most'' irreducible representations. In~\cite{B17}, Bessenrodt proved that all double hooks, i.e. partitions with Durfee size 2, are contained in $\chi^{\rho_k}\otimes \chi^{\rho_k}$, and  Li generalized to triple hooks~\cite{L21}.  Harman and Ryba~\cite{HR22}  showed that $\chi^{\rho_k}\otimes \chi^{\rho_k} \otimes \chi^{\rho_m}$ contains every irreducible representation. 

The semigroup property is a powerful tool to create positive Kronecker triples by induction out of building blocks. However, creating a particular partition, like the staircase or a rectangle, out of building blocks is very restrictive and quickly encounters number theoretic issues. A rectangle can be ``cut'' only into other rectangles. This prompts us to enlarge the class of positive triples. 
In 2020, with Christine Bessenrodt we generalized Conjecture~\ref{conj:tensor_square} as follows (and checked for partitions up to size 20).

\begin{conj}[Bessenrodt-Panova 2020]
For every $n$ there exists a $k(n)$, such that for every $\la \vdash n$ with $\la=\la'$ and $d(\la)>k_n$ which is not the square partition, we have $g(\la,\la,\mu)>0$ for all $\mu \vdash n$.
\end{conj}

Here $d(\la) = \max\{i: \la_i\geq i\}$ is the Durfee square size of the partition. Partial progress on that conjecture will appear in the work of Chenchen Zhao.

Another positivity result, motivated by applications in GCT, was obtained in~\cite{IP} as

\begin{lemma}[{\cite[Thm~1.10]{IP}}]\label{l:IP}
Let \ts $X:= \{1, \. 1^2, \. 1^4, \. 1^6, \. 2\ts 1, \. 3\ts 1\}$,
and let partition \ts $\nu \notin X$.  Denote \ts $\ell := \max\{\ell(\nu)+1,9\}$,
and suppose \ts $r >  3\ell^{3/2}$, \ts $s \geq 3 \ell^2$,
and \ts $|\nu| \leq rs/6$. Let $\nu[N]:= (N-|\nu|,\nu_1,\nu_2,\ldots)$. Then \ts $g(s^r, s^r,\nu[rs])>0$.
\end{lemma}

\medskip

The tensor square conjectures raise the question on simply determining when $g(\la,\mu,\nu)>0$.  It is a consequence of representation theory that when $n>2$ for every $\mu \vdash n$ there is a $\la \vdash n$, such that $g(\la,\la,\mu)>0$, see \cite[Ex. 7.82]{S1} but even that has no combinatorial proof.

To understand why such problems are difficult, we resort  to the computational complexity framework. For Kostka numbers, determining if $K_{\la,\mu}>0$ is easy, it just requires us checking that $\la \succ \mu$ in the dominance order. For Littlewood-Richardson that's a lot less trivial despite the combinatorial rule.
\medskip

\textsf{LRPos}: 

\texttt{Input:} $\la,\mu, \nu$

\texttt{Output:} Is $c^{\la}_{\mu\nu}>0$? 

\medskip

The proof of the Saturation Conjecture by Knutson and Tao showed that an LR coefficient is nonzero if and only if the corresponding hive polytope is nonempty. This polytope is a refinement of the Gelfand-Tsetlin polytope, and is defined by $O(\ell(\la)^2)$ many inequalities. Showing that the polytope is nonempty is thus a linear programming problem, which can be solved in polynomial time, see~\cite{MNS12}. 
\begin{thm}
We have that \textsf{LRPos} $\in \P$ when the input is in binary (and unary). 
\end{thm}

\medskip

\textsf{KronPos}: 

\texttt{Input:} $\la,\mu, \nu$

\texttt{Output:} Is $g(\la,\mu,\nu)>0$? 

\medskip

 In the early stages of GCT Mulmuley conjectured~\cite{MS2} that they would behave like the Littlewood-Richardson, so \textsf{KronPos}$\in \P$, which was recently disproved.

\begin{thm}[\cite{IMW17}]
When the input $\la,\mu,\nu$ is in unary, \textsf{KronPos} is \NP-hard.
\end{thm}

The proof uses the fact that in certain cases $g(\la,\mu,\nu)$ is equal to the number of pyramids with marginals $\la,\mu,\nu$ and deciding if there is such a pyramid is $\NP$-complete. However, the problem is not yet in $\NP$: we do not have, for every $\la,\mu,\nu$, polynomial witnesses showing that $g(\la,\mu,\nu)>0$ when this is the case. This is yet another open problem.

\begin{op}
Determine whether when the input $\la,\mu,\nu$ is in unary, \textsf{KronPos} is in $\NP$. 
\end{op}

Needless to say, the problem would be [computationally] harder when the input is in binary.

\medskip

Positivity of plethysm coefficients is even less understood. Plethysm coefficients $a_\la(d[n])$ are important in GCT and positivity of $a_{\la}(d[n])$ for partitions $\la$ with long first rows was established in~\cite{BIP}. In general, we do not expect the problem to be any easier. 

\medskip

\textsf{PlethPos}: 

\texttt{Input:} $\la,\mu, \nu$

\texttt{Output:} Is $a^\la_{\mu,\nu}>0$? 

\medskip

In~\cite{FI20}, Fischer and Ikenmeyer  showed that \textsf{PlethPos} is $\NP$-hard. As there is no combinatorial interpretation and no nice positivity criteria in general we cannot say whether the problem is in $\NP$.

\begin{op}
Determine whether \textsf{PlethPos} is in $\NP$. 
\end{op}


There is one major conjecture on plethysm coefficients.

\begin{conj}[Foulkes]\label{conj:foulkes}
Let $d > n$, then 
$$a_\la(d[n]) \geq a_\la(n[d])$$
for all $\la \vdash nd$. 
\end{conj}

This conjecture is related to the Alon-Tarsi conjecture, and has appeared in GCT-related research. In~\cite{DIP} we proved it for some families of 3-row partitions.

As most of these quantities can be computed via symmetric group characters, it is natural to ask how easy it is to decide if a character is 0 or not.

\smallskip
\textsf{CharVanish}: 

\texttt{Input}: $\la,\al \vdash n$, unary. 

\texttt{Output}: Is  $\chi^\la(\al)=0?$

\smallskip
\begin{thm}[\cite{IPP22}]\label{thm:char_vanish}
We have that \textsf{CharVanish} is \CeqP-complete under many-to-one reductions.
\end{thm}

The class $\CeqP$ is defined formally as $[ \GapP =0]$, the class of deciding when two $\SP$ functions have equal values, and this is a class believed to be strictly larger than $\coNP$ and $\NP$. Thus we cannot expect that deciding  $\chi=0$ would be in $\coNP$ and $\chi \neq 0$ in $\NP$. This also leads to the proof of Proposition~\ref{prop:char_kron}. See~\cite{Pan23} for a discussion on that topic and sketch of the proof.

\section{Asymptotic Algebraic Combinatorics}

Asymptotic Algebraic Combinatorics is a natural extension of Algebraic Combinatorics inspired by Probability and Asymptotic Representation Theory. Instead of searching for exact formulas, bijections, combinatorial interpretations, we are looking to understand quantities approximately. This is, indeed, what we could hope for with the given structure constants in the absence of nice formulas. 

As mentioned in the Introduction, such problems have natural interplay with Statistical Mechanics and Integrable Probability. They would also be important in the search of multiplicity obstructions in GCT, see~\cite{Pan23}.

\subsection{Maximal multiplicities}

The first step in such asymptotic analysis begins with estimating the range of possible values.
In~\cite{S1,S2}, Stanley observed that 
$$\aligned
\qquad \quad \  \max_{\la\vdash n}  \max_{\mu\vdash n}  \max_{\nu\vdash n}
\,\, g(\la,\mu,\nu) \, & = \, \sqrt{n!} \, e^{-O(\sqrt{n})}\, , \\ 
 \qquad \, \max_{0\le k\le n}  \max_{\la\vdash n}  \max_{\mu\vdash k} 
\max_{\nu\vdash n-k}  \,\, c^\la_{\mu,\nu} \, & = \, 2^{n/2 -  O(\sqrt{n})}.
\endaligned
$$
and asked
\begin{op}[Stanley]
Fix $n$. For which $\la,\mu,\nu$ with $\la \vdash n$, is $c^\la_{\mu\nu}$ asymptotically maximal?
For which $\la,\mu,\nu$ is $g(\la,\mu,\nu)$ asymptotically maximal?
\end{op}

Using various symmetric function and representation theoretic identities involving sums of these coefficients, it was observed in~\cite{PPY1} that such maximal partitions $\la$ must also have asymptotically maximal dimensions $f^\la$. By the pioneering work of Vershik-Kerov~\cite{VK} and Loggan-Shepp~\cite{LoS}, we know that such partitions approach certain explicit curve, which we refer to as the VKLS shape and VKLS partitions, and for all those we have $f^\la = \sqrt{n!} e^{-O(\sqrt{n})}$. 

\begin{thm}[\cite{PPY1}]
Let $\{\la^{(n)}\vdash n\}$, $\{\mu^{(n)}\vdash n\}$, $\{\nu^{(n)}\vdash n\}$ be
three partition sequences, such that
\begin{equation*}
(*) \qquad \qquad g\bigl(\ts \la^{(n)},\ts\mu^{(n)},\ts\nu^{(n)}\ts\bigr) \, = \, \sqrt{n!} \, e^{-O(\sqrt{n})}\..
\end{equation*}
Then $\la^{(n)}, \mu^{(n)}, \nu^{(n)}$ have VKLS shape (i.e. are asymptotically maximal).  Conversely, for every
two VKLS  sequences \ts $\{\la^{(n)}\vdash n\}$ \ts and \ts
$\{\mu^{(n)}\vdash n\}$, there exists a VKLS partition sequence \ts
$\{\nu^{(n)}\vdash n\}$, s.t. (*) holds.
\end{thm}

This statement shows ``existence'', the problem here is to exhibit a particular explicit family of partitions for which this asymptotics holds.

\begin{op}
Determine a family of partitions $\la^{(n)}, \mu^{(n)}, \nu^{(n)} \vdash n$, such that 
$$ g\bigl(\ts \la^{(n)},\ts\mu^{(n)},\ts\nu^{(n)}\ts\bigr) \, = \, \sqrt{n!} \, e^{-O(\sqrt{n})}$$
as $n \to \infty$. 
\end{op}

We believe that for the staircase partitions $\rho_k=(k,k-1,\ldots,1)$ the Kronecker coefficient $g(\rho_k,\rho_k,\rho_k)$ grows superexponentially, and this is one specific candidate which could achieve large, but not maximal, values. However, no exponential lower bounds are known in this case.

\medskip
 
 Similar analysis can be done for the Littlewood-Richardson coefficients. Again, the proofs can only derive ``existence'', but not exhibit any specific families. 

\begin{thm}[\cite{PPY1}] \label{t:LR-intro1}
Fix \ts $0<\theta<1$ \ts and let \ts $k_n := \lfloor\theta n\rfloor$.
Then:
\smallskip

1. \ts for every VKLS partition sequence $\{\la^{(n)}\vdash n\}$,
there exist VKLS partition sequences $\{\mu^{(n)}\vdash k_n\}$
and $\{\nu^{(n)}\vdash n-k_n\}$, s.t.
%
$$
(**) \qquad \ c^{\la^{(n)}}_{\mu^{(n)},\. \nu^{(n)}} \, = \, \binom{n}{k_n}^{1/2} \, e^{-O(\sqrt{n})}\.,
$$
%

2. \ts  for all VKLS partition sequences
$\{\mu^{(n)}\vdash k_n\}$ and $\{\nu^{(n)}\vdash n-k_n\}$, there exists a
VKLS partition sequence $\{\la^{(n)}\vdash n\}$, s.t.~$(**)$ holds,

3. \ts  for all VKLS partition sequences
$\{\la^{(n)}\vdash n\}$ and $\{\mu^{(n)}\vdash k_n\}$, there exists a
partition sequence $\{\nu^{(n)}\vdash n-k_n\}$, s.t.
$$
f^{\nu^{(n)}} \, = \, \sqrt{n!} \, e^{-O(n^{2/3}\log n)} \quad \text{and} \quad
c^{\la^{(n)}}_{\mu^{(n)},\. \nu^{(n)}} \, = \, \binom{n}{k_n}^{1/2} \, e^{-O(n^{2/3}\log n)}\..
$$
\end{thm}

We again ask when such maxima are achieved.

\begin{op}
Let $\theta \in (0,1)$ be some fixed number and let $k_n := \lfloor\theta n\rfloor$.  Find families of partitions $\{\la^{(n)}\vdash n\}$,
 $\{\mu^{(n)}\vdash k_n\}$
and $\{\nu^{(n)}\vdash n-k_n\}$, for which 
$$
  \ c^{\la^{(n)}}_{\mu^{(n)},\. \nu^{(n)}} \, = \, \binom{n}{k_n}^{1/2} \, e^{-O(\sqrt{n})}\.,
$$
\end{op}

While asymptotics of Kronecker coefficients have barely been studied, Schur functions asymptotics in various regimes has been central to the field of Integrable Probability. It is closely related to the Harish-Chandra-Itzhykhson-Zuber integrals in random matrix theory. A step further takes us to Kostka and Littlewood-Richardson coefficients, which can also be estimated using the so-called elliptic integrals, see~\cite{BGH}. There, the problem of asymptotics is translated into an implicit variational analysis, and it still does not give explicit concrete answers. 

\medskip

There are currently no nontrivial asymptotic results of this kind for plethysm coefficients.

\begin{op}
Determine the maximal asymptotic value of the plethysm coefficients 
$$\max_{\la \vdash n^2} a_{\la}(n[n]).$$
\end{op}

The only nontrivial bounds come from their relation to rectangular Kronecker coefficients, and in particular when $\la$ is a two-row partition, see \S~\ref{ss:tight}. We expect that the maxima should be at least exponential. Another source of nontrivial lower bounds can be found in~\cite{FI20} where plethysms are compared to pyramids. Similar comparison was used to obtain a family of large Kronecker coefficients in~\cite{PP20,PP22}. 

\subsection{Bounded rows and diagonals}

In the ``maximal'' regime, the partitions involved were close to the VKLS shape which is symmetric and has length $O(\sqrt{n})$ for partitions of $n$. The asymptotics changes significantly depending on the regimes of convergence for the partitions involved. 

When the lengths of the partitions are bounded, formulas coming from symmetric functions and relating them to contingency tables become useful. We have the following upper bound.

\begin{thm}[\cite{PP20}]
Let $\la,\mu,\nu\vdash n$ such that $\ell(\la)=\ell$, $\ell(\mu) = m$,
and $\ell(\nu)=r$.  Then:
$$g(\la,\mu,\nu) \, \le \, \left(1+\frac{\ell m r}{n}\right)^n
\left(1+\frac{n}{\ell m r}\right)^{\ell m r} .
$$
\end{thm}

In particular,  when $\la=\mu=\nu=( (\ell^2)^\ell) \vdash \ell^3=n$, we have 
$g(\la,\mu,\nu) \leq 4^n$. Yet, a similar lower bound for that case is beyond reach. 

Further analysis using Schur function techniques also gives bounds when the Durfee size (diagonal) of a partition is bounded, even if the length grows. We have the following, see~\cite{PP22}. 

 \begin{thm}[\cite{PP22}]
Let \ts $n, k\ge 1$, and let \.  $\la,\mu,\nu\vdash n$, such that \. $d(\la), \ts d(\mu), \ts d(\nu) \. \le \ts k$.  Then:
\begin{equation}
g(\la,\mu,\nu) \, \leq \, \frac{1}{k^{8k^2}\ts 2^{8k^3}} \,\, n^{4k^3\ts + \ts 13k^2 \ts + \ts 31k}.
\end{equation}
\end{thm}

We also know that the actual maximal value is not much further from this bound. Let
$$
A(n,k) \, := \, \max \. \big\{ \ts g(\la,\mu,\nu) \ : \ \la,\mu,\nu\vdash n \ \
\text{and} \ \ \ell(\la),\ell(\mu),\ell(\nu)\le k \ts \big\}.
$$

\begin{thm}[\cite{PP22}]
\label{t:intro-rows-lower}
For all \ts $k\ge 1$, there is a constant \ts $C_k>0$, such that
\begin{equation}\label{eq:main-lower-A}
A(n,k) \, \geq \, C_k \. n^{k^3\ts - \ts 3k^2 \ts - \ts 3k \ts + \ts 3}  \qquad
\text{for all \ \ $n\ge 1$\ts.}
\end{equation}
\end{thm}

However, we are nowhere close to finding asymptotics close to these maxima for concrete families of partitions.

\begin{op}
For every fixed $k$ determine a family of partitions $\la^{(n)},\mu^{(n)},\nu^{(n)}\vdash n$ with $\ell(\la^{(n)}), \ell(\mu^{(n)}), \ell(\nu^{(n)})\leq k$, such that 
$$g(\la^{(n)},\mu^{(n)},\nu^{(n)}) = O(n^{(k-1)^3}).$$ 
\end{op}

Using the powerful technique of ACSV (Asymptotic Combinatorics in Several Variables), the cases of $k \leq 4$ have been somewhat well understood and explicit asymptotics derived in~\cite{MT}. 
For arbitrary, but fixed, values of $k$ though the ACSV is no longer applicable and we need other ways to compute this.

The same question pertains to the plethysm coefficients.

\begin{op}
Show that  for every fixed $k$ the maximal plethysm coefficient over partitions with length $k$ grows polynomially in $n$. In particular, we have
$$\log \max_{\la \vdash n^2,\, \ell(\la) =k} a_\la(n[n]) = \Omega(\log(n))$$
\end{op}

It was establsihed in~\cite{IP} that plethysms are bounded above by rectangular Kronecker coefficients, 
$$a_\la(d[n]) \leq g(\la,n^d, n^d)$$
for $\ell(\la) \leq m^2$, $\la_1 \geq d(n-m)$ and $d,n$ large enough so we are in the stable regime, i.e. the values of the plethysm and Kronecker are constant when only $d,n$ and $\la_1$ grow. Thus, we would not expect that plethysm coefficients would have larger growth than Kroneckers. 

\subsection{Tight asymptotics}\label{ss:tight}

There are very few nontrivial cases where we have certain explicit formulas for Kronecker and plethysm coefficients and would be able to derive tight asymptotics. One such case is for rectangular Kronecker coefficients $g(\la, m^\ell, m^\ell)$ when $\la$ has two rows. We have that, see~\cite{PPq}, for $\la=(m\ell-n,n)$,
$$
g(m^\ell,m^\ell, \la ) = p_n(\ell,m) \ts - \ts p_{n-1}(\ell,m)  =a_{\la}(m[\ell]),
$$
where $p_n(\ell,m) = \# \{ \al\vdash n, \al_1\leq m, \ell(\al) \leq \ell\}$ is the number of partitions of $k$ whose Young diagram fits inside an $\ell \times m$ rectangle. 

Using this partition counting interpretation, and deriving tight asymptotics for $p_n(m,\ell)$ using tilted geometric random variables, we obtain
\begin{thm}[\cite{MPPe}] 
Let $A := \frac{\ell}{m}$  $B := \frac{n-1}{m^2}$. Let $c, d$ be solutions of a certian system of integral equations with parameters $A,B$.  We have
\[  a_{\la}(m[\ell])= g( \la, m^\ell, m^\ell) = p_{n}(\ell,m) - p_{n-1}(\ell,m) \sim \frac{d}{m}p_{n-1}(\ell,m) \sim \frac{d \,\, e^{m\left[ cA+2dB-\log(1-e^{-c-d})\right]}}{2\pi m^3 \sqrt{D}}. \]
\end{thm} 

Note that some technical conditions for the theorem to hold also include that $\ell,m$ grow at the same rate, and also that $\frac{n}{m\ell}$ should be away from $\frac{1}{2}$. When $n \sim \frac{1}{2} m\ell$ we can no longer estimate the difference between $p_n$'s so tightly and we need a different approach for the following

\begin{op}
Find tight asymptotics when $\la=(\lceil m^2/2  \rceil, \lfloor m^2/2 \rfloor)$ for
$g( m^m, m^m, \la)$.  
\end{op}
The best lower bound in this case comes from characters and is obtained in~\cite{PPq}.

Needless to say, any other explicit tight asymptotic results would be of great interest and, using semigroup property, could lead to improved lower bounds in other regimes.

\section{Reduced Kronecker coefficients}\label{s:reduced}

Since the Kronecker coefficients quickly turned out to be difficult to understand, hopes turned towards the \defn{stable Kronecker coefficients} defined as 

$$\rg(\al,\be,\ga) \. := \. \lim_{n\to \infty} \. g\bigl(\al[n],\be[n],\ga[n]\bigr), \
 \quad \al[n]:= (n-|\al|,\al_1,\al_2,\ldots), \ n\ge |\al|+\al_1.
$$
They generalize the Littlewood-Richardson coefficients, following Theorem~\ref{thm:mur} from~\cite{Mur56}. 
$$
\rg(\al,\be,\ga) \, = \, c^\al_{\be\ga} \quad \text{for} \quad |\al|\. = \. |\be| \ts +\ts |\ga|\ts,
$$
and were thereby called ``extended Littlewood-Richardson coefficients'' in~\cite{Kir}.
They have thus occupied an intermediate spot between the Littlewood-Richardson and ordinary Kronecker coefficients, and have been a subject of independent study and interest, see~\cite{Mur38, Mur56, Bri93, Val99, Kir, Kly, BOR2, BDO, CR, Man15, SS16, IP, PPr, OZ21}. Together with the ordinary Kronecker coefficients, the problem of finding a combinatorial interpretation has also been part of the agenda, see~\cite{Sta00}.

The ordinary Kronecker coefficients can be expressed as a small alternating sum of reduced Kroneckers and reduced Kroneckers are certain sums of ordinary Kroneckers, see~\cite{BOR2,BDO}. These relationships show that reduced Kroneckers are also \SP-hard to compute, see~\cite{PPr}. However, it could not be immediately deduced that their positivity is \NP-hard, which  prompted the following question.

\begin{op}\label{op:reduced}
Given input $\la,\mu,\nu$ (in unary), is deciding whether $\rg(\la,\mu,\nu)>0$ $\NP$-hard? Determine whether computing $\rg(\la,\mu,\nu)$ is \SP-complete? 
\end{op}
As we noted, the second part of that problem would follow if  there is a combinatorial interpretation for them, or for the ordinary Kronecker coefficients. 

Some special cases of combinatorial interpretations can be derived from the ordinary Kronecker coefficients. Separately, in~\cite{CR} a combinatorial interpretation was given when $\mu,\nu$ are rectangles and $\la$ is one row, and other similar cases were derived in~\cite{BO}.  Methods to compute them, as well as for the ordinary Kronecker coefficients, have been developed in a series of papers, see~\cite{BDO,BOR2,OZ20,OZ21}. 

Until recently, it was believed they are strictly better behaved than the ordinary Kroneckers, with the following conjecture appearing in~\cite{Kir,Kly}

\begin{conj}[Kirillov, Klyachko]\label{conj:sat}
The reduced Kronecker coefficients satisfy the saturation property:
$$\rg(N\al,N\be, N\ga) > 0 \quad \text{for some} \. \ N\ge 1
\quad \Longrightarrow \quad  \rg(\al,\be,\ga) > 0\ts.
$$
\end{conj}

Using the recent positivity results, in particular for rectangular Kroneckers from~\cite{IP}, this was also disproved.
\begin{thm}[\cite{PPr}]
For all \ts $k\ge 3$, the triple of partitions
\ts $\bigl(1^{k^2-1},1^{k^2-1},k^{k-1}\bigr)$ \ts is a counterexample
to the Conjecture.  For every partition~$\ga$ \ts s.t.\ \ts
$\ga_2\ge 3$,
there are infinitely many pairs \ts $(a,b)\in \mathbb{N}^2$ \ts s.t.
\ts $(a^b,a^b,\ga)$ \ts is a counterexample.
\end{thm}
As shown in~\cite{PPr}, asymptotically, and complexity-wise they behave more similarly to the ordinary Kronecker than to the Littlewood-Richardson coefficients.

\begin{rem}
In the course of writing this paper, the first part of Open Problem~\ref{op:reduced} was  settled in~\cite{IP23}, where it was shown that every Kronecker coefficient is equal to a [not much larger] reduced Kronecker coefficient:
$$g(\la,\mu,\nu) = \rg( \nu_1^{\ell(\la)}+\la, \nu_1^{\ell(\mu)}+\mu, (\nu_1^{\ell(\la)+\ell(\mu)},\nu) ).$$
 Thus the reduced Kronecker coefficients are just as hard as the Kroneckers, and the complexity results follow from\S~\ref{s:combin}. This also shows that the problems are essentially the same, and the reduced Kronecker coefficients are not easier in general. 
\end{rem}

\section{Methods}

While the concrete open problems inspire work in the area, the pertinent task is to create tools for understanding the Kronecker and plethysm coefficients. Of course, an ultimate tool could be the ``combinatorial interpretation'', but this itself has become the hardest problem to understand and its solution requires new tools. So far, as very little progress has been made in this direction, it is clear that we need something new. 

Kronecker coefficients can be interpreted via $S_n$ as multiplicities. Via Schur-Weyl duality, they are also dimensions of invariant spaces under $GL$ action. The last interpretation is what gives the semigroup property, something that cannot be seen from the $S_n$ perspective. On the other hand, the transpositional invariance $g(\la,\mu,\nu) = g(\la',\mu',\nu)$  is visible from their $S_n$ definition. There are many more examples showing the variety of results one can get using different methods and the limitations of each one also. Perhaps the best tool would be able to unify all these views. 

Symmetric functions have been one of the most useful tools in Algebraic Combinatorics, and is certainly the author's favorite. The Littlewood-Richardson coefficients are structure constants in the ring of symmetric functions and thus have been very successfully studied using Schur functions.  But the Kronecker coefficients cannot be interpreted in that way. Motivated by that, a family of symmetric functions was developed in~\cite{OZ20, OZ21}. Symmetric functions have not been very useful in establishing positivity results, because the formulas involved are often alternating sums. However, they are useful for finding asymptotics and bounds as done in~\cite{PP20,PP22} for example. Within integrable probability they have been a lot more successful, from Asymptotic Representation Theory, through lozenge tilings, to Random Matrices, see e.g.~\cite{BG15, GP, BGH} and references therein. 

Plethysm coefficients can also be studied using symmetric functions, yet again positivity properties cannot usually be derived this way. For example in~\cite{DIP} symmetric functions were used to derive some explicit formulas for $a_{\la}(d[n])$ with $\la$ a three row partition, and later this method generalized in \cite{FI20} to exhibit an alternating formula and put them in $\GapP$. 

A technique to show positivity is the construction of explicit highest weight tableau evaluations, see~\cite{Ful00}, which is easier for plethysms but also applicable for Kroneckers. For example, this technique was used in~\cite{I15} to prove Theorem~\ref{thm:rho_mu}, generalized in~\cite{LS}, and in~\cite{BIP} to show positive plethysm coefficients. 

Useful tools in Representation Theory are induction and restriction of modules, passing to the representation theory of other groups. For example, the character bound~\ref{prop:char_kron} was obtained by restricting to the alternating group $A_n \subset S_n$, see~\cite{BB,PPV, PPq}. Similar methods were developed to extract positivity criteria working with representations over finite fields, see e.g.~\cite{BBL}. In another direction, diagram algebras can also be used to study Kronecker, see e.g.~\cite{BDO}, and plethysm coefficients, see e.g.~\cite{COSSZ,OSSZ}.

Another useful approach comes from discrete geometry. This is not at all surprising given that many objects in representation theory correspond to integer points in polytopes, e.g. SSYTs are integer points in Gelfand-Tsetlin polytopes, see e.g.~\cite{DM06} for some complexity connections also. Vallejo used such relations to develop methods of discrete tomography and find a lower bound for Kronecker coefficients via pyramids, see~\cite{Val97,Val99,IMW17}. Kronecker coefficients are bounded above by 3d contingency tables with given 2d marginals, as seen in~\cite{IMW,PP20} to obtain upper bounds. Similar expressions tying them to Kostant partition functions were developed in~\cite{MRS21} and used in~\cite{MT} to find asymptotics. 

Finally, Computational Complexity theory provides some tools besides the framework. The framework is to explain the hardness of the problems and what possible answers to expect. The tools come in the form of gadgets and basic hard/complete problems to reduce to like \textsf{3SAT}, \textsf{PerfectMatching}, \textsf{SetPartition} etc. Embedding such a complete problem in our problem is an art, very similar to finding injections, see e.g.~\cite{GJ} for some classical results. Other tools come from algebraic complexity and GCT, where the complexity of certain polynomials could imply inequalities between the relevant structure constants.  

\medskip
   
Despite the variety of questions, methods and applications, the Kronecker and plethysm coefficients are still a mysterious black box. It sometimes gives us presents. But we cannot see inside... just yet.


\begin{thebibliography}{Bou12300}

\bibitem[Aar16]{Aa}
S.~Aaronson, $\P\overset{?}=\NP$, in \emph{Open problems in mathematics},
Springer, Cham, 2016, 1--122.






\bibitem[BO05]{BO} C.~Ballantine, R.~Orellana,  A combinatorial interpretation for the coefficients in the Kronecker product $s_{(n-p,p)}*s_\la$. \emph{S\'em. Lothar. Combin.} {\bf 54A} (2005/07), 29 pp.




\bibitem[BGH]{BGH} S.~Belinschi, A.~Guionnet, J.~Huang,
Large deviation principles via spherical integrals,
\emph{Prob. Math. Phys.}{\bf 3} (3) (2022), pp. 543--625.

\bibitem[BB04]{BB} C.~Bessenrodt and C.~Behns,
On the Durfee size of Kronecker products of characters of
the symmetric group and its double covers,
{\em J.~Algebra}~{\bf 280} (2004), 132--144.



\bibitem[BB17]{BB17} 
C.~Bessenrodt, C.~Bowman, Multiplicity-free Kronecker products of characters of the symmetric groups. \emph{Adv. Math.} {\bf 322} (2017), pp. 473--529. 

\bibitem[B18]{B17} C.~Bessenrodt, Critical classes, Kronecker products of spin characters, and the Saxl conjecture, \emph{Alg. Combinatorics} {\bf 1} (2018), 353--369.


\bibitem[BBL]{BBL} C. Bessenrodt, C. Bowman, L. Sutton, Kronecker positivity and 2-modular representation theory, \emph{Trans. Amer. Math. Soc. Ser. B}{\bf 8}(2021), 1024--1055.



\bibitem[BBR06]{BBR06}
F.~Bergeron,  R.~Biagioli and M.~H.~Rosas,
Inequalities between Littlewood--Richardson coefficients,
\emph{J.~Combin.\ Theory, Ser.~A }~\textbf{113} (2006), 567--590.


\bibitem[Bia98]{Bia}
Ph.~Biane, Representations of Symmetric Groups and
Free Probability, \emph{Adv.\ Math.}~\textbf{138} (1998), 126--181.



\bibitem[B13]{B13}  J.~Blasiak, Kronecker coefficients for one hook shape,\emph{ S\'em. Lothar. Combin}, {\bf 77}:2016--2017.

\bibitem[BL18]{BL} J.~Blasiak, R.I.~Liu, Kronecker coefficients and noncommutative super Schur functions. \emph{ J. Combin. Theory Ser. A}{\bf 158} (2018), 315--361.

\bibitem[BMS15]{BMS}
J.~Blasiak, K.~Mulmuley, M.~Sohoni, Geometric complexity theory IV: nonstandard quantum group for the Kronecker problem. \emph{Mem. Amer. Math. Soc.}{\bf 235} (2015), no. 1109.


\bibitem[BI18]{BI}
M.~Bl\"aser and C.~Ikenmeyer,
\emph{Introduction to geometric complexity theory},
Summer school lecture notes, 2018, 148~pp.;
{\tt https://tinyurl.com/nhe2wxvw}

\bibitem[BDO15]{BDO} C.~Bowman, M.~De~Visscher, R.~Orellana, The partition algebra and the Kronecker coefficients. \emph{Trans. Amer. Math. Soc.} {\bf 367} (2015), no. 5, 3647--3667. 

\bibitem[BOR11]{BOR2}
E.~Briand, R.~Orellana and M.~Rosas,
The stability of the Kronecker product of Schur functions,
\emph{J.\ Algebra}~\textbf{331} (2011), 11--27.

\bibitem[Bri93]{Bri93} M.~Brion,
Stable properties of plethysm: on two conjectures of Foulkes,
\emph{Manuscripta mathematica} (1993),
{\bf 80}:347–371.

\bibitem[BG15]{BG15} A.~Bufetov, V.~Gorin,  Representations of classical Lie groups and quantized free convolution. \emph{Geom. Funct. Anal.} {\bf 25}, 763--814 (2015). 

\bibitem[B\"{u}r00a]{Bur00} P.~B\"urgisser,
\emph{Completeness and reduction in algebraic complexity theory},
Springer, Berlin, 2000, 168~pp.

\bibitem[B\"{u}r00b]{Bur0} P.~B\"urgisser, Cook's versus Valiant's
hypothesis, \emph{Theor.\ Comp.\ Sci.}~{\bf 235} (2000), 71--88.



\bibitem[BCS97]{BCS}
P.~B\"urgisser, M.~Clausen and M.~A.~Shokrollahi,
\emph{Algebraic complexity theory}, Springer, Berlin, 1997.

\bibitem[BI08]{BI08} P.~Bürgisser, C.~Ikenmeyer,  The complexity of computing Kronecker coefficients. \emph{ 20th Annual International Conference on Formal Power Series and Algebraic Combinatorics} (FPSAC 2008), pp. 357--368.



\bibitem[BIP19]{BIP}
P.~B\"urgisser, C.~Ikenmeyer and G.~Panova,
No occurrence obstructions in geometric complexity theory,
\emph{J.~Amer.~Math.~Society}~\textbf{32} (2019), 163--193; extended abstract in
\emph{57-th FOCS} (2016), IEEE, Los Alamitos, CA, 386--5.

\bibitem[BLMW11]{BLMW}
P.~B\"urgisser, J.~M.~Landsberg, L.~Manivel and J.~Weyman,
An overview of mathematical issues arising in the Geometric Complexity Theory
approach to $\VP\overset{?}=\VNP$, \emph{SIAM J.\ Comput.} {\bf 40} (2011), 1179--1209.

\bibitem[CHM07]{CHM}
 M.~Christandl, A.~Harrow, G.~Mitchison,  Nonzero Kronecker coefficients and what they tell us about spectra. \emph{Comm. Math. Phys.} {\bf 270 }(2007), no. 3, pp.575--585. 

\bibitem[CDW12]{CDW}
M.~Christandl, B.~Doran and M.~Walter, Computing Multiplicities of
Lie Group Representations, in \emph{Proc.\ 53-rd FOCS} (2012), IEEE, 639--648.

\bibitem[COSSZ22]{COSSZ}
L.~Colmenarejo, R.~Orellana, F.~Saliola, .~Schilling, M.~Zabrocki,
The mystery of plethysm coefficients, \url{https://arxiv.org/abs/2208.07258} (2022).

\bibitem[CR15]{CR}
L.~Colmenarejo and M.~Rosas,
Combinatorics on a family of reduced Kronecker coefficients,
\emph{C.R.\ Math.\ Acad.\ Sci.\ Paris}~\textbf{353} (2015), no.~10, 865--869.





\bibitem[DM06]{DM06}
J.~A.~De~Loera and T.~B.~McAllister,
On the computation of Clebsch--Gordan coefficients and the dilation effect,
\emph{Experimental Math.}~\textbf{15} (2006), 7--19.



\bibitem[D93]{Dvir}
Y.~Dvir,  On the Kronecker product of $S_n$ characters,
\emph{J.~Algebra}~\textbf{154} (1993), 125--140.



\bibitem[DIP19]{DIP}
J.~D\"orfler, C.~Ikenmeyer and G.~Panova,
{On geometric complexity theory: Multiplicity obstructions are stronger than occurrence obstructions} in
\emph{Proc.\ 46-th ICALP} (2019), Art.~51, 14~pp.
i




\bibitem[FI20]{FI20}
N.~Fischer and C.~Ikenmeyer,
The computational complexity of plethysm coefficients,
\emph{Comput.\ Complexity}~\textbf{29} (2020), no.~2,
Paper~8, 43~pp.

\bibitem[FRT54]{FRT}
J.~S.~Frame, G.~de~B. Robinson and R.~M. Thrall,
The hook graphs of the symmetric group,
{\em Canad.\ J.\ Math.}~\textbf{6} (1954), 316--324.


\bibitem[Ful97]{Ful}
W.~Fulton, \emph{Young tableaux},
Cambridge Univ.\ Press, Cambridge, UK, 1997, 260~pp.

\bibitem[Ful00]{Ful00}
W.~Fulton, Eigenvalues, invariant factors, highest weights, and Schubert calculus,
\emph{Bull.\ AMS}~\textbf{37} (2000), 209--249.



\bibitem[GJ79]{GJ}
M.~R.~Garey and D.~S.~Johnson, \emph{Computers and intractability},
Freeman, San Francisco, CA, 1979.

\bibitem[GR85]{GR85} A.~M.~Garsia and J.~Remmel, Shuffles of permutations and the Kronecker product, \emph{Graphs and
Combinatorics} (1985), 1:217-263.

\bibitem[GIP17]{GIP}
 F.~Gesmundo, C.~Ikenmeyer and G.~Panova,
{Geometric complexity theory and matrix powering},  \emph{Diff.\ Geom.\ Applications}~{\bf 55} (2017), 106--127.




\bibitem[GP15]{GP}
V.~Gorin and G.~Panova,
{Asymptotics of symmetric polynomials with applications to
statistical mechanics and representation theory},
\emph{Ann.\ Probab.}~{\bf 43} (2015), 3052--3132.


\bibitem[HR22]{HR22} N.~Harman, C.~Ryba, A Tensor-Cube Version of the Saxl Conjecture, \texttt{arXiv:2206.13769 }. 

\bibitem[HSTZ13]{HSTZ}
G.~Heide, J.~Saxl, P.~H.~Tiep, A.~Zalesski,
Conjugacy action, induced representations and the Steinberg square for simple groups of Lie type, 
\emph{Proc. Lond. Math. Soc.} (3) {\bf 106} (2013), no. 4, 908--930.

\bibitem[I15]{I15}
C.~Ikenmeyer, The Saxl conjecture and the dominance order, \emph{Disc. Math.} (11) {\bf 6} (2015), pp.1970--1975.

\bibitem[IMW17]{IMW17}
C.~Ikenmeyer, K.~D.~Mulmuley and M.~Walter,
On vanishing of Kronecker coefficients,
\emph{Computational Complexity}~\textbf{26} (2017), 949--992.


\bibitem[IP22]{IP22}
C.~Ikenmeyer and I.~Pak,
What is in~$\SP$ and what is not?,  preprint (2022),
82~pp.; extended abstract to appear in \emph{Proc.\ 63rd FOCS} (2022); \ts {\tt arXiv:2204.13149}.

\bibitem[IPP22]{IPP22}
C.~Ikenmeyer, I.~Pak and G.~Panova,
Positivity of the symmetric group characters is as hard as the polynomial time hierarchy,
preprint (2022), 15 pp.; extended abstract to appear in \emph{Proc.\ 34th SODA} (2023); \ts {\tt arXiv:2207.05423}.

\bibitem[IP17]{IP}
C.~Ikenmeyer and G.~Panova, Rectangular Kronecker coefficients and plethysms
in geometric complexity theory, \emph{Adv.\ Math.}~\textbf{319} (2017), 40--66;
extended abstract in \emph{57th FOCS} (2016), IEEE, Los Alamitos, CA, 396--4055.

\bibitem[IP23]{IP23}
C.~Ikenemeyer and G.~Panova, All Kronecker coefficients are reduced Kronecker coefficients, (2023). 


%

\bibitem[Kir04]{Kir}
A.~N.~Kirillov,
An invitation to the generalized saturation conjecture,
\emph{Publ.\ RIMS}~\textbf{40} (2004), 1147--1239.

\bibitem[Kly04]{Kly}
A.~Klyachko,
Quantum marginal problem and representations of the symmetric group,
{\tt arXiv:quant-ph/0409113}.

\bibitem[KT99]{KT99}
A.~Knutson and T.~Tao,
The honeycomb model of $\GL_n(\cc)$ tensor products~I:
Proof of the saturation conjecture,
{\em J.~AMS} {\bf 12} (1999), 1055--1090.


\bibitem[Las79]{Las79} A. Lascoux, Produit de Kronecker des repr\'esentations du groupe sym\'etrique. \emph{ S\'eminaire d\`Alg\`ebre Paul Dubreil et
Marie-Paule Malliavin: Proceedings}, Paris 1979 (32\`eme Ann\'ee), pages 319–329 (1979).

\bibitem[L21]{L21} Xin~Li, Saxl Conjecture for triple hooks, \emph{Discrete Math.}
{\bf 344} (2021) 6.

\bibitem[Liu17]{Liu17} Ricky Liu,
A simplified Kronecker rule for one hook shape,
\emph{Proc.  Amer. Math.
Soc.}(2017), {\bf 145}(9):3657-3664.

\bibitem[LoS]{LoS}
B.~F.~Logan and L.~A.~Shepp,
A variational problem for random Young tableaux,
\emph{Adv.\ Math.}~\textbf{26} (1977), 206--222.

\bibitem[LS17]{LS}
S.~Luo,M.~Sellke, The Saxl conjecture for fourth powers via the semigroup property. \emph{J. Algebraic Combin.} {\bf 45} (2017), no. 1, 33--80.



\bibitem[Mac95]{Mac}
I.~G.~Macdonald,
\emph{Symmetric functions and Hall polynomials} (Second edition),
Oxford University Press, New York, 1995.

\bibitem[Man15]{Man15} L~Manivel,
On the asymptotics of Kronecker coefficients,
\emph{J.  Alg. Combinatorics}(2015),
{\bf 42}(4):999-1025.




\bibitem[MPP19a]{MPPe}
S.~Melczer, G.~Panova and R.~Pemantle,
{Counting partitions inside a rectangle},
\emph{SIAM J.\ Discrete Math.}~\textbf{34} (2020), 2388--2410;
extended abstact in {\em Proc.\ 31st FPSAC}, 2019.






\bibitem[MRS21]{MRS21} M.~Mishna, M.~Rosas, S.~Sundaram, Vector partition functions and Kronecker coefficients,
\emph{ J. Phys. A: Math. Theor.} (2021) {\bf 54}.

\bibitem[MT]{MT} M.~Mishna, S.~Trandafir, Estimating and computing Kronecker Coefficients: a vector partition function approach, \texttt{arXiv:2210.12128}. 






\bibitem[Mul11]{Mul11}
K.~D.~Mulmuley, Geometric Complexity Theory VI: The flip via positivity,
preprint (2011), 40~pp., \ts available at \ts
{\tt http://gct.cs.uchicago.edu/gct6.pdf};
cf.\ {\tt arXiv:0704.0229}, 139~pp.


\bibitem[Mul17]{Mul} K.~Mulmuley,
Geometric Complexity Theory V.\ Efficient algorithms for Noether normalization,
\emph{J.\ AMS}~{\bf 30} (2017), 225--309.

\bibitem[MNS12]{MNS12}
K.~D.~Mulmuley, H.~Narayanan and M.~Sohoni,
Geometric complexity theory~III.  On deciding nonvanishing of a
Littlewood-Richardson coefficient, \emph{J.\ Algebraic Combin.}~\textbf{36}
(2012), 103--110.

\bibitem[MS01]{MS1}
K.~D.~Mulmuley and M.~Sohoni, Geometric complexity theory.~{I} \ts
{A}n approach to the {\P} vs.\ {\NP}  and related problems,
{\em SIAM J.\ Comput.}~\textbf{31}, 2001, 496--526.

\bibitem[MS08]{MS2}
K.~D.~Mulmuley and M.~Sohoni, Geometric complexity theory.~{II} \ts
{T}owards explicit obstructions for embeddings among class varieties,
{\em SIAM J.\ Comput.}~\textbf{38}, 2008, 1175--1206.


\bibitem[Mur38]{Mur38}
F.~D.~Murnaghan,
The analysis of the Kronecker product of irreducible representations of the symmetric group,
\emph{Amer.\ J.\ Math.}~\textbf{60} (1938), 761--784.

\bibitem[Mur56]{Mur56}
F.~D.~Murnaghan,
On the Kronecker product of irreducible representations of the symmetric group,
\emph{Proc.\ Natl.\ Acad.\ Sci.\ USA}~\textbf{42} (1956), 95--98.

\bibitem[Nar06]{Nar}
H.~Narayanan,
On the complexity of computing Kostka numbers and Littlewood--Richardson coefficients,
\emph{J.~Algebraic Combin.}~\textbf{24} (2006), 347--354.





\bibitem[OSSZ]{OSSZ}
R.~Orellana, F.~Saliola, A.~Schilling, M.~Zabrocki,  Plethysm and the algebra of uniform block permutations, \emph{Alg. Combinatorics} {\bf 5} (2022) (5), pp. 1165--1203.

\bibitem[OZ20]{OZ20} R.~Orellana, M.~Zabrocki,
 A combinatorial model for the decomposition of multivariate polynomials rings as an $S_n$-module
\emph{ The Elect. Journal of Comb.}{\bf 27}(3) (2020), \#P3.24.


\bibitem[OZ21]{OZ21} R.~Orellana, M.~Zabrocki,
Symmetric group characters as symmetric functions
\emph{Adv. Math.} {\bf 390} (2021). 


\bibitem[Pak22+]{Pak22}
I.~Pak, What is a combinatorial interpretation?, preprint (2022), 58 pp.;
to appear in \emph{Proc.\ OPAC}, AMS, Providence, RI.

\bibitem[PP13]{PP13}
I.~Pak and G.~Panova, Strict unimodality of $q$-binomial coefficients,
\emph{C.\ts\/R.\/ Math.\/ Acad.\/ Sci.\/ Paris} \textbf{351} (2013), 415--418.

\bibitem[PP14]{PP14}
I.~Pak and G.~Panova,
Unimodality via Kronecker products, \emph{J.~Algebraic Combin.}~\textbf{40}
(2014), 1103--1120.

\bibitem[PP17a]{PPq} I.~Pak, G.~Panova, Bounds on certain classes of Kronecker and q-binomial coefficients,
\emph{Journal of Combinatorial Theory, Series A}{\bf 147} (2017), pp.1--17.

\bibitem[PP17b]{PP17}
I.~Pak and G.~Panova,
On the complexity of computing Kronecker coefficients,
\emph{Comput.\ Complexity}~\textbf{26} (2017), 1--36.

\bibitem[PP20a]{PP20}
I.~Pak and G.~Panova,  Bounds on Kronecker coefficients via contingency tables,
\emph{Linear Alg.\  Appl.}~\textbf{602} (2020), 157--178.

\bibitem[PP20b]{PPr} I.~Pak, G.~Panova, Breaking down the reduced Kronecker coefficients,
\emph{Comptes Rendus. Math\'ematique} {\bf 358} (4), 463-468

\bibitem[PP23]{PP22} I.~Pak, G.~Panova, Durfee squares, symmetric partitions and bounds on Kronecker coefficients, \emph{J. Algebra} (2023), to appear. 

\bibitem[PPV16]{PPV} I.~Pak, G.~Panova, E.~Vallejo, Kronecker products, characters, partitions, and the tensor square conjectures,
\emph{Advances in Mathematics} {\bf 288} (2016), pp.702--731.

\bibitem[PPY19a]{PPY1}
I.~Pak, G.~Panova and D.~Yeliussizov,
On the largest Kronecker and Littlewood--Richardson coefficients,
\emph{J.~Combin.\ Theory, Ser.~A}~\textbf{165} (2019), 44--77.





\bibitem[Pan15]{Pporto} G.~Panova, Kronecker coefficients: combinatorics, complexity and beyond, AMS-EMS meeting, Porto, Portugal (2015), slides at \url{https://tinyurl.com/bdd4w3pj}


\bibitem[Pan23]{Pan23} G.~Panova, Computational Complexity in Algebraic Combinatorics, \emph{Current Developments in Mathematics}, Harvard University (2023).   



\bibitem[RW94]{RW} J.~Remmel, T.~Whitehead,  On the Kronecker product of Schur functions of two row shapes. \emph{Bull. Belg. Math. Soc. Simon Stevin} {\bf 1}(1994), no. 5, 649--683.

\bibitem[Sag01]{Sag}
B.~E.~Sagan, \emph{The symmetric group} (Second ed.), Springer, New York, 2001.

\bibitem[SS16]{SS16} S.~V.~Sam and A.~Snowden, \emph{Proof of Stembridge’s conjecture on stability of Kronecker coefficients},
\emph{J. Alg. Combinatorics}(2016), {\bf 43}(1):1-10.




\bibitem[Sta99]{S1}
R.~P.~Stanley, {\em Enumerative Combinatorics}, vol.~1 (Second ed.)
and~2, Cambridge Univ.\ Press, 2012 and~1999.

\bibitem[Sta00]{Sta00}
R.~P.~Stanley,
Positivity problems and conjectures in algebraic combinatorics,
in \emph{Mathematics: frontiers and perspectivies},
AMS, Providence, RI, 2000, 295--319.

\bibitem[Sta20]{S2} R.~P.~Stanley, Supplementary Excercies to \cite{S1}, available at \url{http://www-math.mit.edu/~rstan/ec}.

\bibitem[T15]{T15} V.~Tewari, Kronecker coefficients for some near-rectangular partitions, \emph{J. Algebra} (2015) {\bf 429}, p. 287--317.


\bibitem[V79a]{V1}
L.~G.~Valiant, Completeness classes in algebra. \emph{Proc.\ 11th  STOC} (1979), 249--261.

\bibitem[V79b]{V2}
L.~G.~Valiant,   The  complexity  of  computing  the  permanent.
\emph{Theor.\  Comp.\ Sci.}~{\bf 8} (1979), 189--201.



\bibitem[Val97]{Val97}
E.~Vallejo, Reductions of additive sets, sets of uniqueness and pyramids,
\emph{Discrete Math.}~\textbf{173} (1997), 257--267.

\bibitem[Val99]{Val99}
E.~Vallejo,
Stability of Kronecker products of irreducible characters of
the symmetric group,
\emph{Electron.\ J.\ Combin.}~\textbf{6} (1999), RP~39, 7~pp.

\bibitem[VK]{VK}
A.~M.~Vershik and S.~V.~Kerov,
Asymptotic of the largest and the typical dimensions of irreducible
representations of a symmetric group, \emph{Funct.\ Anal.\ Appl.}~\textbf{19} (1985), 21--31.

\bibitem[Wig19]{Wig}
A.~Wigderson, \emph{Mathematics and Computation}, monograph draft, 2019.




\end{thebibliography}
\end{document}